\documentclass[11pt, oneside, a4paper]{article}
\usepackage{cite}
\usepackage{amsfonts}
\usepackage{fancyhdr,amsthm,amsmath,amssymb,mathrsfs}
\newtheorem{thm}{\bf Theorem}[section]
\newtheorem{lem}[thm]{Lemma}
\newtheorem{prop}[thm]{Proposition}
\newtheorem{cor}[thm]{Corollary}

\theoremstyle{definition}
\newtheorem{defn}[thm]{Definition}
\newtheorem{rem}[thm]{Remark}

\numberwithin{equation}{section}
\newenvironment{prf}{\noindent {\bf Proof.}\rm}{\qed}
\newcommand{\pg}{{\rm Path}(G)}

\newcommand{\cir}[1]{{\rm Cir}(#1)}

\newcommand{\cnw}{{\rm C}(W)}
\newcommand{\cnh}[1]{{\rm C}(#1)}

\newcommand{\In}{{\rm Inv}(G)}

\newcommand{\al}{\alpha}
\newcommand{\be}{\beta}
\newcommand{\ga}{\gamma}
\newcommand{\ze}{\zeta}

\begin{document}
\title{\Large\bf Congruences on graph inverse semigroups\thanks{Partially
supported by the Fundamental Research Funds for the Central
Universities (XDJK2016B038).}}
\author{Zheng-Pan Wang \\
{\small\em School of Mathematics and Statistics, Southwest University, $400715$ Chongqing, China} \\
{\small\em E-mail: zpwang@swu.edu.cn} }
\date{}
\maketitle

\begin{abstract}
Inverse graph semigroups were defined by Ash and Hall in 1975. They found necessary and sufficient conditions for the semigroups to be congruence free. In this paper we give a description of congruences on a graph inverse semigroup in terms of the underlying graph. As a consequence, we show that the inverse graph semigroup of a finite graph is congruence Noetherian.
\medskip


\end{abstract}

\section{Introduction and main results}

In 1975, graph inverse semigroups were first introduced by Ash and Hall in \cite{AshHall} where they gave necessary and sufficient conditions for the semigroups to be congruence free. Such semigroups are defined based on special kinds of directed graphs where there exists at most one edge from one vertex to another. The definition of graph inverse semigroups based on general directed graphs can be found in \cite{Jones, JonesLawson, MesyanMitchell}. These inverse semigroups are related to various types of algebras, for instance, $C^\ast$-algebras, Cohn path algebras, Leavitt path algebras and so on (see \cite{AbramsPino, AlaliGilbert, JonesLawson, Krieger, KumjianPaskRaeburnRenault, KumjianPaskRaeburn, Paterson}). Graph inverse semigroups have also been studied in their own right in recent years (see \cite{AlaliGilbert, ChihPlessas, Jones, JonesLawson, MesyanMitchell, MesyanMitchellMoraynePeresse}).

A {\it directed graph} $G = (V, E, r, s)$ consists of $V, E$ and functions $r, s : E \rightarrow V$. The elements of $V$ are called {\it vertices} and the elements of $E$ are called {\it edges}. For each edge $e$, $s(e)$ is the {\it source} of $e$ and $r(e)$ is the {\it range} of $e$. The cardinality $|s^{-1}(v)|$ of a vertex $v$ is called the {\it index} of $v$. In this paper, we shall refer to a directed graphs simply as a ``graph''.

Now let $G = (V, E, r, s)$ be a graph. The {\it graph inverse semigroup} $\In$ of $G$ is the semigroup with zero generated by $V$ and $E$, together with the set $E^\ast$ of variables $\{e^\ast : e \in E\}$, satisfying the following relations for all $u, v \in V$ and $e, f \in E$:
\begin{description}
\item[{\rm(GI1)}] $u v = \delta_{u,v} u$;
\item[{\rm(GI2)}] $s(e) e = e r(e) = e$;
\item[{\rm(GI3)}] $r(e) e^\ast = e^\ast s(e) = e^\ast$;
\item[{\rm(CK1)}] $f^\ast e = \delta_{f,e} r(e)$.
\end{description}
A {\it path} in $G$ is a sequence $\al = e_1 \cdots e_n$ of edges such that $r(e_i) = s(e_{i+1})$ for $i = 1, \cdots, n-1$. In such a case, $s(\al) := s(e_1)$ is the {\it source} of $\al$, $r(\al) := r(e_n)$ is the {\it range} of $\al$, $|\al| := n$ is the {\it length} of $\al$ and $\al^\ast := e_n^\ast \cdots e_1^\ast$ the {\it ghost path} of $\al$. We view an element $v$ of $V$ as a path of length $0$, $s(v) = r(v) = v$ and $v^\ast = v$. Denote as $\pg$ the set of all paths in $G$. Every nonzero element in $\In$ can be uniquely written as $\al\be^\ast$ for some $\al, \be \in \pg$ with $r(\al) = r(\be)$. One can also verify that $\In$ is an inverse semigroup with zero such that $(\al\be^\ast)^{-1} = \be\al^\ast$ for any $\al, \be \in \pg$. The idempotent set $E(\In)$, which is a (meet) semilattice, consists of all elements $\al\al^\ast$ with $\al \in \pg$. Moreover, for any $\al \be^\ast, \ze \eta^\ast \in \In$, we have
\begin{equation} \label{multiplication}
(\al \be^\ast) (\ze \eta^\ast) =
\left\{ \begin{array}{ll}
(\al\xi)\eta^\ast & \mbox{if $\ze = \be\xi$ for some $\xi \in \pg$} \\
\al(\eta\xi)^\ast & \mbox{if $\be = \ze\xi$ for some $\xi \in \pg$} \\
0 & \mbox{otherwise}.
\end{array}\right.
\end{equation}

In \cite{MesyanMitchell}, Mesyan and Mitchell showed that the quotient of any graph inverse semigroup by a Rees congruence is always isomorphic to another graph inverse semigroup. They also found a graph theoretical condition for a not necessarily Rees congruence and a graph theoretical condition for a graph inverse semigroup to have only Rees congruences. Jones also investigated non-Rees congruences on certain kinds of graph inverse semigroups (see \cite{Jones}).

In this paper we give a description of congruences on a graph inverse semigroup in terms of the underlying graph that is different from \cite{MesyanMitchell}.  This description implies that the inverse graph semigroup of a finite graph is congruence Noetherian.

To present our main result, we need some more notions and notation. {\it In what follows, $G = (V, E, r, s)$ is always a graph and $\In$ is the graph inverse semigroup of $G$}.

A subset $H$ of $V$ is {\it hereditary} if $s(e) \in H$ always implies $r(e) \in H$ for $e \in E$. Let $H$ be a hereditary subset of $V$, $V_1 = V \setminus H$, $E_1 = \{e \in E : r(e) \not\in H\}$, $s_{G \setminus H} = s |_{E_1}$ and $r_{G \setminus H} = r |_{E_1}$. Then we denote the graph $(V_1, E_1, s_{G \setminus H}, r_{G \setminus H})$ as $G \setminus H$.

Let $\rho$ be a congruence on $\In$. Then clearly $H = 0 \rho \cap V$ is a hereditary subset of $V$. Passing to the graph $G \setminus H$, we shall assume without loss of generality that $H = \emptyset$. So for graph inverse semigroups, we only need to describe these {\it special congruences} where no vertex is equivalent to 0.

Let $\al = e_1 \cdots e_n$ be a path in $G$ and $v(\al)$ denote the set $\{s(e_i) : i = 1, \cdots, n\}$. If $| \al | = 0$, we appoint that $v(\al) = \emptyset$. An edge $e$ is an {\it exit} to the path $\al$ if there exists $i$ such that $s(e) = s(e_i)$ and $e \neq e_i$. The path $\al$ is said to be {\it no-exit} if it does not have an exit. The path $\al$ is called a {\it cycle} if $s(\al) = r(\al)$ and $s(e_i) \neq s(e_j)$ whenever $i \neq j$. For any cycle $c$, let $\cir c$ be the set of all cycles that are obtained from $c$ by cyclic permutations. Two cycles $c_1$ and $c_2$ are said to be {\it disjoint} if $v(c_1) \cap v(c_2) = \emptyset$.

We write $\bar{V} := \{v \in V : | s^{-1}(v) | = 1\}$, the set of vertices of index one. For any $v \in \bar{V}$, we denote the unique edge in $s^{-1}(v)$ as $e_v$. Let $W$ be a subset of $\bar{V}$, $\mathbb{Z}^+$ be the set of all positive integers and $\cnw$ be the set of all cycles whose vertices lie in $W$. Since all vertices in $W$ have index one, an arbitrary cycle $c \in \cnw$ is no-exit and any two cycles in $\cnw$ are either disjoint or cyclic permutations of each other.

\begin{defn} \label{Pair}
A {\it cycle function} $f : \cnw \rightarrow \mathbb{Z}^+ \cup \{\infty\}$ is a function that is invariant under cyclic permutations. A {\it congruence pair $(W, f)$ of $G$} consists of a subset $W$ of $\bar{V}$ and a cycle function $f$.
\end{defn}

\begin{defn} \label{Triple}
Let $H$ be a hereditary subset of $V$ and $(W, f)$ be a congruence pair of $G \setminus H$. Then $(H, W, f)$ is called a {\it congruence triple of $G$}.
\end{defn}

Let us show how a congruence gives rise to a congruence triple (respectively, a special congruence gives rise to a congruence pair).

Let $\rho$ be a congruence on $\In$, $H = 0 \rho \cap V$ and $W = \{v \in V \setminus H : e  e^\ast \rho v \mbox{ for } e \in s^{-1}_{G \setminus H} (v)\}$. Later we shall see from Lemma~\ref{NoExit} that in this case all vertices in $W$ have index one in $G \setminus H$. Let $\cnw$ consist of all distinct cycles whose vertices lie in $W$. For $c \in \cnw$, let $f(c)$ be the smallest positive integer $m$ such that $c^m \rho s(c)$. If no power of $c$ is equivalent to $s(c)$, then we define $f(c) = \infty$. Then $T(\rho) = (H, W, f)$ is a congruence triple.

If $\rho$ is a special congruence, then $T(\rho) = (\emptyset, W, f)$. In this case, $(W, f)$ is a congruence pair.

Conversely, let $(W, f)$ be a congruence pair of $G$ and let $\wp(W, f)$ denote the congruence generated by $\mathbf{R}$ consisting of all pairs $(e_v e_v^\ast, v)$ for $v \in W$ and $(c^{f(c)}, s(c))$ for $c \in \cnw$ with $f(c) \in \mathbb{Z}^+$.

Let $(H, W, f)$ be a congruence triple of $G$ and let $\wp(H, W, f)$ denote the congruence generated by all pairs $(v, 0)$ for $v \in H$, $(ee^\ast, w)$ for $w \in W$ with $s(e) = w$ and $(c^{f(c)}, s(c))$ for $c \in \cnw$ with $f(c) \in \mathbb{Z}^+$.

\begin{thm} \label{Main1}
The mapping $T$ from the set of all special congruences on $\In$ to the set of all congruence pairs of $G$ and the mapping $\wp$ from the set of all congruence pairs of $G$ to the set of all special congruences on $\In$ are inverses. In particular, there exists a one-to-one correspondence between special congruences on $\In$ and congruence pairs of $G$.
\end{thm}

Theorem~\ref{Main1} immediately implies

\begin{thm} \label{Main2}
The mapping $T$ from the set of all congruences on $\In$ to the set of all congruence triples of $G$ and the mapping $\wp$ from the set of all congruence triples of $G$ to the set of all congruences on $\In$ are inverses. In particular, there exists a one-to-one correspondence between congruences on $\In$ and congruence triples of $G$.
\end{thm}

A graph inverse semigroup is said to be {\it congruence Noetherian} if for any congruence sequence
$$
\rho_1 \subseteq \rho_2 \subseteq \cdots \subseteq \rho_n \subseteq \cdots,
$$
there exists $m \in \mathbb{Z}^+$ such that $\rho_m = \rho_{m+k}$ for all $k \in \mathbb{Z}^+$.

\begin{cor} \label{Noetherian}
Each graph inverse semigroup of a finite graph is congruence Noetherian.
\end{cor}

\section{Proof of the main results}

We need to prove Theorem~\ref{Main1} and Corollary~\ref{Noetherian}. Let $\mathbb{N}$ be the set of all nonnegative integers. For any cycle $c$ and any positive integer $m$, we appoint that $c^0 = s(c)$ and $(c^\ast)^m = c^{-m}$. The next lemma directly follows from the related definitions.

\begin{lem} \label{CircleAdd}
For any no-exit cycle $c$ and $\al \in \pg$ with $s(c) = s(\al)$, there exist $\al_1, \al_2 \in \pg$ and $k \in \mathbb{N}$ such that $\al = c^k \al_1$, $c = \al_1 \al_2$ and $|\al_1| < |c|$.
\end{lem}

\begin{lem} \label{Circle}
For any no-exit cycle $c$, $\al \in \pg$ and $k \in \mathbb{Z}^+$, if $s(c) = s(\al)$, then there exists $c_1 \in \cir c$ such that
\begin{description}
\item[{\rm(1)}] $\al^\ast c^k \al = c_1^k$, and
\item[{\rm(2)}] $c^k \al \al^\ast = \al c_1^k \al^\ast$.
\end{description}
\end{lem}

\begin{prf}
We see from Lemma~\ref{CircleAdd} that $\al = c^l \al_1$ where $l \in \mathbb{N}$ and $\al_1 \in \pg$ such that $c = \al_1 \al_2$ for some $\al_2 \in \pg$. Obviously, $c_1 = \al_2 \al_1 \in \cir c$. It follows from (CK1) that
\begin{align*}
& \al^\ast c^k \al = \al_1^\ast (c^\ast)^l c^k c^l \al_1 = \al_1^\ast c^k \al_1 = \al_2 c^{k-1} \al_1 = c_1^k, \mbox{ and} \\
& c^k \al \al^\ast = c^k c^l \al_1 \al^\ast = c^l \al_1 c_1^k \al^\ast = \al c_1^k \al^\ast.
\end{align*}
That is, (1) and (2) hold.
\end{prf}
\medskip

Now let $\rho$ be a special congruence on $\In$ and $T(\rho) = (W, f)$. We shall show that $\rho = \wp (W, f)$.
\bigskip

The next Lemma directly follows from (GI1), (CK1) and the definition of a congruence on a semigroup.

\begin{lem} \label{Basic}
Let $u, v \in V$ and $\al, \be, \ga \in \pg$.
\begin{description}
\item[{\rm(1)}] If $u \neq v$, then $(u, v) \not\in \rho$;
\item[{\rm(2)}] if $(\al \be^\ast, v) \in \rho$, then $s(\al) = s(\be) = v$;
\item[{\rm(3)}] $(\al \be, s(\al)) \in \rho$ if and only if $(\be \al, s(\be)) \in \rho$;
\item[{\rm(4)}] if $(\al \be \ga^\ast \al^\ast, s(\al)) \in \rho$, then $(\be \ga^\ast, s(\be) \in \rho$;
\item[{\rm(5)}] if $(\al \al^\ast, s(\al)) \in \rho$, $(\be \ga^\ast, s(\be)) \in \rho$ and $r(\al) = s(\be)$, then $(\al \be \ga^\ast \al^\ast, s(\al)) \in \rho$.
\end{description}
\end{lem}

\begin{lem} \label{NoExit}
For any $\al \be^\ast \in \In$, if $(\al \be^\ast, s(\al)) \in \rho$, then both $\al$ and $\be$ are no-exit.
\end{lem}

\begin{prf}
Suppose that $\al = e_1 \cdots e_n$ and $e$ is an exit to $\al$ with $s(e) = s(e_i)$. It follows from Lemma~\ref{Basic}(2) and (CK1) that
$$
(e_i \cdots e_n \be^\ast e_1 \cdots e_{i-1}, s(e_i)) \in \rho.
$$
We obtain from (GI3) and (CK1) that
$$
(0, r(e)) = (0, e^\ast e) = (e^\ast e_i \cdots e_n \be^\ast e_1 \cdots e_{i-1} e, e^\ast s(e_i) e) \in \rho
$$
which contradicts to the hypothesis $0\rho = \{0\}$. Hence, $\al$ is no-exit. We see that $\be$ is also no-exit since $(\al \be^\ast, s(\al)) \in \rho$ implies that $(\be \al^\ast, s(\be) \in \rho$.
\end{prf}

\begin{lem} \label{Subpath}
For any  $\al \in \pg$, if $(\al \al^\ast, s(\al)) \in \rho$ and $\al = \xi \be \eta$ for some $\xi, \be, \eta \in \pg$, then $(\be \be^\ast, s(\be)) \in \rho$.
\end{lem}

\begin{prf}
We see from the hypothesis that $(\xi \be \eta \eta^\ast \be^\ast \xi^\ast, s(\xi)) \in \rho$. Then it follows from  Lemma~\ref{Basic}(4) that $(\be \eta \eta^\ast \be^\ast, s(\be)) \in \rho$ and $(\eta \eta^\ast, s(\eta)) \in \rho$. We see from $(\eta \eta^\ast, s(\eta)) \in \rho$ and the compatibility of $\rho$ that $(\be \eta \eta^\ast \be^\ast, \be \be^\ast) \in \rho$. Hence, we obtain by the transitivity of $\rho$ that $(\be \be^\ast, s(\be)) \in \rho$.
\end{prf}  \medskip

For any $c \in \cnw$, the semigroup $\langle c, c^\ast \rangle$ generated by $c$ and $c^\ast$ is the (infinite) bicyclic semigroup (see \cite{Howie}), and $\langle c, c^\ast \rangle / \rho$ is either the bicyclic semigroup or a cyclic group of order $f(c)$. This implies

\begin{lem} \label{Divisible}
For any $c \in \cnw$, if $f(c) \in \mathbb{Z}^+$, then for every $m \in \mathbb{Z}^+$, $(c^m, s(c)) \in \rho$ implies that $f(c) \mid m$.
\end{lem}

A {\it closed path based at $v$} is a path $\al = e_1 \cdots e_n$ such that $s(\al) = r(\al) = v$. A closed path $\al =  e_1 \cdots e_n$ based at $v$ is called a {\it closed simple path based at $v$} if $s(e_i) \neq v$ for every $i > 1$. A cycle is a closed simple path.

\begin{lem}[{\cite[Lemma~2.3]{AbramsPino}}] \label{CSP}
For any $v \in V$ and  any closed path $\al$ based at $v$, there exist unique closed simple paths $\eta_1, \cdots, \eta_m$ based at $v$ such that $\al = \eta_1 \cdots \eta_m$.
\end{lem}

\begin{lem} \label{PowerofCycle}
If $\al \in \pg \setminus V$ and $(\al, s(\al)) \in \rho$, then there exists a no-exit cycle $c$ such that $\al = c^m$ for some positive integer $m$.
\end{lem}

\begin{prf}
We see from Lemma~\ref{Basic}(2) that $\al$ is a closed path based at $s(\al)$. It follows from Lemma~\ref{CSP} that there exist unique closed simple paths $\mu_1, \cdots, \mu_m$ based at $s(\al)$ such that $\al = \mu_1 \cdots \mu_m$. Then we obtain from Lemma~\ref{NoExit} that $\mu_1 = \cdots = \mu_m = c$ for some no-exit cycle $c$.
\end{prf}  \medskip

The following lemma follows from \cite[Proposition~8]{MesyanMitchell}.

\begin{lem} \label{Classofv}
For any $v \in V$, every element in $v \rho$ is of the form $\al \eta \al^\ast$ or $\al \eta^\ast \al^\ast$ for some $\al, \eta \in \pg$ such that $s(\al) = v$ and $r(\al) = s(\eta) = r(\eta)$.
\end{lem}

\begin{lem} \label{Elementofvrho}
For any $v \in V$, every element in $v \rho$ is of the form either $\al \al^\ast$ or $\al c^{mf(c)} \al^\ast$, where $\al \in \pg$, $s(\al) = v$, $v(\al) \subseteq W$ and in the second case, $m \in \mathbb{Z} \setminus \{0\}$, $r(\al) = s(c)$, $c \in \cnw$ with $f(c) \in \mathbb{Z}^+$.
\end{lem}

\begin{prf}
Note Lemma~\ref{Classofv} and take $x$ in $v \rho$. Then either $x = \al \eta \al^\ast$ or $x = \al \eta^\ast \al^\ast$. We already have $r(\al) = s(\eta)$, $s(\al) = v$ and $v(\al) \subseteq W$ by Lemma~\ref{NoExit}. If $|\eta| = 0$, then it is clear $x = \al \al^\ast$. If $|\eta| > 0$, then noticing $\eta = \eta (s(\eta))^\ast$, we see from Lemmas~\ref{Basic}(4), (2), \ref{PowerofCycle} and \ref{Divisible} that $x = \al c^{m f(c)} \al^\ast$ where $m \in \mathbb{Z} \setminus \{0\}$ and $c \in \cnw$ with $f(c) \in \mathbb{Z}^+$.
\end{prf}

\begin{prop} \label{CPC}
$\wp(W, f) = \rho$.
\end{prop}

\begin{prf}
Recall that $\mathbf{R}$ consists of all pairs $(e_v e_v^\ast, v)$ for $v \in W$ and $(c^{f(c)}, s(c))$ for $c \in \cnw$ with $f(c) \in \mathbb{Z}^+$. Clearly, $\mathbf{R} \subseteq \rho$ so that $\wp(W, f) \subseteq \rho$. For any $\al \be^\ast, \xi \eta^\ast \in \In \setminus \{0\}$, if $(\al \be^\ast, \xi \eta^\ast) \in \rho$, we may suppose $|\al| \leq |\xi|$ without loss of generality. Then we see from Lemma~\ref{NoExit} that $\xi = \al \xi_1$ for some $\xi_1 \in \pg$. Concerning $\be$ and $\eta$, we have the following two cases. \medskip \\
\textbf{Case 1.} $\eta = \be \eta_1$ for some $\eta_1 \in \pg$. It follows from $(\al \be^\ast, \xi \eta^\ast) \in \rho$ and (CK1) that $(r(\al), \xi_1 \eta_1^\ast) \in \rho$. Then we see from Lemmas~\ref{Elementofvrho} and \ref{Subpath} that $(r(\al), \xi_1 \eta_1^\ast) \in \wp(W, f)$. Hence, we have $(\al\be^\ast, \xi \eta^\ast) = (\al \be^\ast, \al \xi_1 \eta_1^\ast \be^\ast) \in \wp(W, f)$. \medskip \\
\textbf{Case 2.} $\be = \eta \be_1$ for some $\be_1 \in \pg \setminus V$. It again follows from $(\al \be^\ast, \xi \eta^\ast) \in \rho$ and (CK1) that $(r(\al), \xi_1 \be_1) \in \rho$. Then we observe from Lemma~\ref{PowerofCycle} that there exists $c \in \cnw$ with $f(c) \in \mathbb{Z}^+$ and such that $\xi_1 \be_1 = c^m$ for some $m \in \mathbb{Z}^+$. Thus, it follows from Lemma~\ref{Divisible} that there exists $k \in \mathbb{Z}^+$ such that $m = k f(c)$. Now we get from the definition of $\wp(W, f)$ that $(c^{kf(c)}, s(c)) \in \wp(W, f)$ which means that $(r(\al), \xi_1 \be_1) \in \wp(W, f)$. Noticing (CK1), we have
$$
(\xi_1^\ast, \be_1) \in \wp(W, f) \Longrightarrow (\be_1^\ast, \xi_1) \in \wp(W, f) \Longrightarrow (\al \be_1^\ast \eta^\ast, \al \xi_1 \eta^\ast) \in \wp(W, f).
$$
That is, $(\al \be^\ast, \xi \eta^\ast) \in \wp(W, f)$. \medskip \\
As a consequence, we obtain that $\rho \subseteq \wp(W, f)$ which leads to $\rho = \wp(W, f)$.
\end{prf}  \medskip

Next let $(W, f)$ be a congruence pair of $G$ and $\wp(W, f) = \varrho$. We shall prove that $\varrho$ is a special congruence and $T(\varrho) = (W, f)$.
\bigskip

Let $S$ be a semigroup, let $\mathbf{T}$ be a relation on $S$ and denote the congruence generated by $\mathbf{T}$ as $\mathbf{T}^\sharp$. If $c, d$ in $S$ are such that $c = xay, d = xby$ for some $x, y \in S^1$, where either $(a, b) \in \mathbf{T}$ or $(b, a) \in \mathbf{T}$, we say that $c$ is connected to $d$ by an {\it elementary $\mathbf{T}$-transition} and write $c \rightarrow d$.

\begin{rem} \label{Trans}
Recall that the relation $\mathbf{R}$ consists of all pairs $(e_v e_v^\ast, v)$ for $v \in W$ and $(c^{f(c)}, s(c))$ for $c \in \cnw$ with $f(c) \in \mathbb{Z}^+$ and that $\wp(W, f) = \varrho = \mathbf{R}^\sharp$. The definition of an elementary $\mathbf{R}$-transition $c \rightarrow d$ in $\In$ can be adjusted as $c = xay, d = xby$ for some $x, y \in \In$ (not adjoined with $1$) where either $(a, b) \in \mathbf{R}$ or $(b, a) \in \mathbf(R)$. In fact, if either $x$ or $y$ equals to $1$, then either $x$ or $y$ can be replaced by some vertex since we have (GI1--3).
\end{rem}

\begin{lem} \label{Transition}
Let $\mathbf{T}$ be a relation on a semigroup $S$ and $a, b \in S$. Then $(a, b) \in \mathbf{T}^\sharp$ if and only if either $a = b$ or for some $n \in \mathbb{Z}^+$, there exists a sequence
$$
a = z_1 \rightarrow z_2 \rightarrow \cdots \rightarrow z_n = b
$$
of elementary $\mathbf{T}$-transitions connecting $a$ to $b$.
\end{lem}

\begin{lem} \label{vinWsigma}
For any $v \in V$, every element in $v \varrho$ is of the form either $\al \al^\ast$ or $\al c^{kf(c)} \al^\ast$, where $\al \in \pg$, $s(\al) = v$, $v(\al) \subseteq W$, and in the second case, $k \in \mathbb{Z}$, $r(\al) = s(c)$, $c \in \cnw$ with $f(c) \in \mathbb{Z}^+$.
\end{lem}

\begin{prf}
If $w \in v \varrho$, then we know from Lemma~\ref{Transition} that there exists a sequence
$$
v = z_1 \rightarrow z_2 \rightarrow \cdots \rightarrow z_n = w.
$$
of elementary $\mathbf{R}$-transitions. We now prove by induction on $n$. Obviously, the result is true if $n = 1$. Suppose that $z_{n-1}$ has one of the two forms in the lemma and that $z_{n-1} = \al \be^\ast x \xi \eta^\ast, z_n = \al \be^\ast y \xi \eta^\ast$ for some $\al \be^\ast, \xi \eta^\ast \in \In \setminus \{0\}$, where either $(x, y) \in \mathbf{R}$ or $(y, x) \in \mathbf{R}$. Then we have the following cases and subcases. \medskip \\
\textbf{Case 1.} $z_{n-1} = \ga \ga^\ast$ for some $\ga \in \pg$. Note $(\ref{multiplication})$ and that $v(\ga) \subseteq W$. \\
\textbf{Subcase 1.1.} $x \in W$ and $y = e_x e_x^\ast$. \\
\textbf{Subcase 1.1.1.} $\be = \xi \be_1$ for some $\be_1 \in \pg$. Then we get $z_{n-1} = \al \be_1^\ast \eta^\ast = \ga \ga^\ast$ which means that $\al = \ga = \eta \be_1$. Thus, if $|\be| = 0$, then we see that $z_n = \ga e_x e_x^\ast \ga^\ast$. Clearly, $v(\ga e_x) \subseteq W$. If $|\be| > 0$, then we know that $z_n = \al \be^\ast e_x e_x^\ast \xi \eta^\ast = \al \be^\ast \xi \eta^\ast = z_{n-1}$.  \\
\textbf{Subcase 1.1.2.} $\xi = \be \xi_1$ for some $\xi_1 \in \pg$. Then we get $z_{n-1} = \al \xi_1 \eta^\ast = \ga \ga^\ast$ which leads to $\al \xi_1 = \ga = \eta$. One can similarly obtain that either $z_n = \ga e_x e_x^\ast \ga^\ast$ or $z_n = z_{n-1}$ via discussing the two possibilities $|\xi| = 0$ and $|\xi| > 0$. \\
\textbf{Subcase 1.2.} $x \in W$ and $y = c^{f(c)}$ for some $c \in \cnw$ with $s(c) = x$ and $f(c) \in \mathbb{Z}^+$.  \\
\textbf{Subcase 1.2.1.} $\be = \xi \be_1$ for some $\be_1 \in \pg$. As in Subcase~1.1.1, we get $\al = \ga = \eta \be_1$. It follows from Lemma~\ref{Circle}(1) that $z_n = \al \be_1^\ast \xi^\ast c^{f(c)} \xi \eta^\ast = \al \be_1^\ast c_1^{f(c_1)} \eta^\ast$ for some $c_1 \in \cir c$. By Lemma~\ref{CircleAdd}, there exist $\be_2, \be_3 \in \pg$ and $l \in \mathbb{N}$ such that $\be_1 = c_1^l \be_2$ and $c_1 = \be_2 \be_3$. So we get $z_n = \eta c_1^l \be_2 \be_2^\ast (c_1^\ast)^l c_1^{f(c_1)} \eta^\ast$. Take $c_2 = \be_3 \be_2$. Clearly, $c_2 \in \cir c$. If $l \geq f(c_1) = f(c_2)$, then we have
$$
z_n = \eta c_1^l \be_2 \be_2^\ast (c_1^\ast)^{l - f(c_1)} \eta^\ast = \eta \be_2 c_2^l (c_2^\ast)^{l - f(c_2)} \be_2^\ast \eta^\ast = (\eta \be_2 c_2^{l - f(c_2)}) c_2^{f(c_2)} ((c_2^\ast)^{l - f(c_2)} \be_2^\ast \eta^\ast).
$$
Furthermore, $v(\eta \be_2 c_2^{l - f(c_2)}) \subseteq W$. If $l < f(c_1)$, then we have
$$
z_n = \eta c_1^l \be_2 \be_2^\ast c_1^{f(c_1) - l} \eta^\ast = \eta c_1^l \be_2 \be_2^\ast \be_2 \be_3 (c_1)^{f(c_1) - l - 1} \eta^\ast = \eta c_1^{f(c_1)} \eta^\ast.
$$
Obviously, $v(\eta) \subseteq W$. \\
\textbf{Subcase 1.2.2.} $\xi = \be \xi_1$ for some $\xi_1 \in \pg$. As in Case~1.1.2, we also get $\al \xi_1 = \ga = \eta$. It follows from Lemma~\ref{Circle}(1) that $z_n = \al \be^\ast c^{f(c)} \be \xi_1 \eta^\ast = \al c_1^{f(c_1)} \xi_1 \eta^\ast$ for some $c_1 \in \cir c$. So we see from Lemma~\ref{Circle}(2) that
$$
z_n = \al c_1^{f(c_1)} \xi_1 \xi_1^\ast \al^\ast = \al \xi_1 c_2^{f(c_2)} \xi_1^\ast \al^\ast
$$
for some $c_2 \in \cir c$. We already have $v(\al \xi_1) \subseteq W$. \\
\textbf{Subcase 1.3.} $x = e_y e_y^\ast$ and $y \in W$. Then we get $z_{n-1} = \al \be^\ast e_y e_y^\ast \xi \eta^\ast$. If $|\be| \neq 0$ or $|\xi| \neq 0$, then we easily see that $z_{n-1} = \al \be^\ast \xi \eta^\ast = z_n$. So we only need to consider the case where $|\be| = |\xi| = 0$. In this case, we have $\al e_y e_y^\ast \eta^\ast = \ga \ga^\ast$. So we get $\al e_y = \ga = \eta e_y$ which means that $z_n = \al \al^\ast$ with $v(\al) \subseteq W$. \\
\textbf{Subcase 1.4.} $x = c^{f(c)}$ and $y \in W$ for some $c \in \cnw$ with $s(c) = y$ and $f(c) \in \mathbb{Z}^+$.  \\
\textbf{Subcase 1.4.1.} $\be = \xi \be_1$ for some $\be_1 \in \pg$. Then we get $z_{n-1} = \al \be_1^\ast \xi^\ast c^{f(c)} \xi \eta^\ast$ and $z_n = \al \be_1^\ast \eta^\ast$. It follows from Lemma~\ref{Circle}(1) that there exists $c_1 \in \cir c$ such that $z_{n-1} = \al \be_1^\ast c_1^{f(c_1)} \eta^\ast$. By Lemma~\ref{CircleAdd}, we see that $\be_1 = c_1^l \be_2$ where $l \in \mathbb{N}$ and $c_1 = \be_2 \be_3$ for some $\be_3 \in \pg$. Clearly, $c_2 = \be_3 \be_2 \in \cir c$. If $l \geq f(c)$, then we obtain that $z_{n-1} = \al \be_2^\ast (c_1^\ast)^{l - f(c)} \eta^\ast = \ga \ga^\ast$ which leads to $\al = \ga = \eta c_1^{l - f(c)} \be_2$. Hence, we have
$$
z_n = \eta c_1^{l - f(c_1)} \be_2 \be_2^\ast (c_1^\ast)^l \eta^\ast = \eta \be_2 c_2^{l - f(c_2)} (c_2^\ast)^l \be_2^\ast \eta^\ast = \ga (c_2^\ast)^{f(c_2)} \ga^\ast
$$
as required. If $l < f(c_1)$, then we see that $z_{n-1} = \al \be_2^\ast c_1^{f(c_1) - l} \eta^\ast = \al \be_3 c_1^{f(c_1) - l -1} \eta^\ast = \ga \ga^\ast$ which leads to $\al \be_3 c_1^{f(c_1) - l -1} = \ga = \eta$. Hence, we have
$$
z_n = \al \be_2^\ast (c_1^\ast)^l (c_1^\ast)^{f(c_1) - l - 1} \be_3^\ast \al^\ast = \al (c_2^\ast)^{f(c_2)} \al^\ast.
$$
Obviously, $v(\al) \subseteq W$.  \\
\textbf{Subcase 1.4.2.} $\xi = \be \xi_1$ for some $\xi_1 \in \pg$. Then we get $z_{n-1} = \al \be^\ast c^{f(c)} \be \xi_1 \eta^\ast$ and $z_n = \al \xi_1 \eta^\ast$. It follows from Lemma~\ref{Circle}(1) that there exists $c_1 \in \cir c$ such that $z_{n-1} = \al c_1^{f(c_1)} \xi_1 \eta^\ast$. So we see that $\al c_1^{f(c)} \xi_1 = \ga = \eta$. Therefore we obtain from Lemma~\ref{Circle}(2) that
$$
z_n = \al \xi_1 \xi_1^\ast (c_1^\ast)^{f(c_1)} \al^\ast = \al (c_1^{f(c_1)} \xi_1 \xi_1^\ast)^\ast \al^\ast = \al \xi_1 (c_2^\ast)^{f(c_2)} \xi_1^\ast \al^\ast
$$
for some $c_2 \in \cir c$. Moreover, $v(\al \xi_1) \subseteq W$. \medskip \\
\textbf{Case 2.} $z_{n-1} = \ga c^{kf(c)} \ga^\ast$, where $\ga \in \pg$, $s(\ga) = v$, $v(\ga) \subseteq W$, $k \in \mathbb{Z}^+$, $r(\ga) = s(c)$ and $c \in \cnw$ with $f(c) \in \mathbb{Z}^+$.
\\
\textbf{Subcase 2.1.} $x \in W$ and $y = e_x e_x^\ast$. \\
\textbf{Subcase 2.1.1.} $\be = \xi \be_1$ for some $\be_1 \in \pg$. Then we get $z_{n-1} = \al \be_1^\ast \eta^\ast = \ga c^{kf(c)} \ga^\ast$ which means that $\al = \ga c^{kf(c)}$ and $\ga = \eta \be_1$. Thus, if $|\be| = 0$, then we see from Lemma~\ref{Circle}(2) that $z_n = \ga c^{kf(c)} e_x e_x^\ast \ga^\ast = \ga e_x c_1^{kf(c_1)} e_x^\ast \ga^\ast$ for some $c_1 \in \cir c$. Obviously, $v(\ga e_x) \subseteq W$. If $|\be| > 0$, then we know that $z_n = \al \be^\ast e_x e_x^\ast \xi \eta^\ast = \al \be^\ast \xi \eta^\ast = z_{n-1}$.  \\
\textbf{Subcase 2.1.2.} $\xi = \be \xi_1$ for some $\xi_1 \in \pg$. Then we get $z_{n-1} = \al \xi_1 \eta^\ast = \ga c^{kf(c)} \ga^\ast$ and similarly obtain that either $z_n = \ga e_x c_1^{kf(c_1)} e_x^\ast \ga^\ast$ for some $c_1 \in \cir c$ or $z_n = z_{n-1}$ via discuss the two possibilities $|\xi| = 0$ and $|\xi| > 0$. \\
\textbf{Subcase 2.2.} $x \in W$ and $y = d^{f(d)}$ for some $d \in \cnw$ with $s(d) = x$ and $f(d) \in \mathbb{Z}^+$. It is not difficult to see that $c$ and $d$ are cyclic permutations of each other since they are no-exit.  \\
\textbf{Subcase 2.2.1.} $\be = \xi \be_1$ for some $\be_1 \in \pg$. As in Subcase~2.1.1, we get $\al = \ga c^{kf(c)}$ and $\ga = \eta \be_1$. It follows from Lemma~\ref{Circle}(1) that $z_n = \al \be_1^\ast \xi^\ast d^{f(d)} \xi \eta^\ast = \al \be_1^\ast c_1^{f(c_1)} \eta^\ast$ for some $c_1 \in \cir c$. By Lemma~\ref{CircleAdd}, there exist $\be_2, \be_3 \in \pg$ and $l \in \mathbb{N}$ such that $\be_1 = c_1^l \be_2$ and $c_1 = \be_2 \be_3$. Moreover, we get $z_n = \ga c^{kf(c)} \be_2^\ast (c_1^\ast)^l c_1^{f(c_1)} \eta^\ast$. So we see that $c = \be_3 \be_2$ since $r(\be_2) = s(c)$. Hence, we see that
$$
z_n = \eta c_1^l \be_2 c^{kf(c)} \be_2^\ast (c_1^\ast)^l c_1^{f(c_1)} \eta^\ast = \eta \be_2 c^{kf(c)+l} \be_2^\ast (c_1^\ast)^l c_1^{f(c_1)} \eta^\ast.
$$
If $l \geq f(c)$, then we have
\begin{align*}
z_n & = \eta \be_2 c^{kf(c)+l} \be_2^\ast (c_1^\ast)^{l-f(c)} \eta^\ast = \eta \be_2 c^{kf(c)+l} (c^\ast)^{l-f(c)} \be_2^\ast \eta^\ast \\
& = (\eta \be_2 c^{l - f(c)}) c^{(k+1) f(c)} ((c^\ast)^{l - f(c)} \be_2^\ast \eta^\ast).
\end{align*}
We observe that $v(\eta \be_2 c_1^{l - f(c)}) \subseteq W$. If $l < f(c)$, then we have
$$
z_n = \eta c_1^{kf(c_1)+l} \be_2 \be_2^\ast c_1^{f(c_1)-l} \eta^\ast = \eta c_1^{kf(c_1)+l} \be_2 \be_2^\ast \be_2 \be_3 c_1^{f(c_1) - l - 1} \eta^\ast = \eta c_1^{(k+1)f(c_1)} \eta^\ast.
$$
Obviously, $v(\eta) \subseteq W$. \\
\textbf{Subcase 2.2.2.} $\xi = \be \xi_1$ for some $\xi_1 \in \pg$. As in Subcase~2.1.2, we also get $z_{n-1} = \al \xi_1 \eta^\ast = \ga c^{kf(c)} \ga^\ast$ so that $\al \xi_1 = \ga c^{kf(c)}$ and $\eta = \ga$. It follows from Lemma~\ref{Circle}(1) that $z_n = \al \be^\ast d^{f(d)} \be \xi_1 \eta^\ast = \al c_1^{f(c_1)} \xi_1 \eta^\ast$ for some $c_1 \in \cir c$. By Lemma~\ref{CircleAdd}, there exist $\xi_2, \xi_3 \in \pg$ and $l \in \mathbb{N}$ such that $\xi_1 = c_1^l \xi_2$ and $c_1 = \xi_2 \xi_3$. Let $c_2 = \xi_3 \xi_2$. Clearly, $c_2 \in \cir {c_1}$ and $\xi_1 = \xi_2 c_2^l$. It follows from $\al \xi_1 = \ga c^{kf(c)}$ that $c = c_2$. Thus we observe that
$$
z_n = \al c_1^{f(c_1)+l} \xi_2 \eta^\ast = \al \xi_2 c^{f(c)+l} \eta^\ast = \al \xi_1 c^{f(c)} \eta^\ast = \ga c^{(k+1)f(c)} \ga^\ast.
$$
We already have $v(\ga) \subseteq W$. \\
\textbf{Subcase 2.3.} $x = e_y e_y^\ast$ and $y \in W$. Then we get $z_{n-1} = \al \be^\ast e_y e_y^\ast \xi \eta^\ast$. If $|\be| \neq 0$ or $|\xi| \neq 0$, then we easily see that $z_{n-1} = \al \be^\ast \xi \eta^\ast = z_n$. So we only need to consider the case where $|\be| = |\xi| = 0$ which means that $\al e_y e_y^\ast \eta^\ast = \ga c^{kf(c)} \ga^\ast$. So we get $\al e_y = \ga c^{kf(c)}$ and $\ga = \eta e_y$ which lead to $\al e_y = \eta e_y c^{kf(c)} = \eta c_1^{kf(c_1)} e_y$. Hence we get $\al = \eta c_1^{kf(c)}$ so that $z_n = \al \eta^\ast = \eta c_1^{kf(c_1)} \eta^\ast$ with $v(\eta) \subseteq W$. \\
\textbf{Subcase 2.4.} $x = d^{f(c)}$ and $y \in W$ for some $d \in \cnw$ with $s(d) = y$ and $f(d) \in \mathbb{Z}^+$. Similar to Subcase~2.2, we observe that $c$ and $d$ are cyclic permutations of each other.  \\
\textbf{Subcase 2.4.1.} $\be = \xi \be_1$ for some $\be_1 \in \pg$. Then we get $z_{n-1} = \al \be_1^\ast \xi^\ast d^{f(d)} \xi \eta^\ast$ and $z_n = \al \be_1^\ast \eta^\ast$. It follows from Lemma~\ref{Circle}(1) that there exists $c_1 \in \cir c$ such that $z_{n-1} = \al \be_1^\ast c_1^{f(c_1)} \eta^\ast$. By Lemma~\ref{CircleAdd}, we see that $\be_1 = c_1^l \be_2$ where $l \in \mathbb{N}$ and $c_1 = \be_2 \be_3$ for some $\be_3 \in \pg$. If $l \geq f(c_1)$, then we obtain that $z_{n-1} = \al \be_2^\ast (c_1^\ast)^{l - f(c_1)} \eta^\ast = \ga c^{kf(c)} \ga^\ast$ which leads to $\al = \ga c^{kf(c)}$ and $\ga = \eta c_1^{l - f(c_1)} \be_2$. These also imply that $c = \be_3 \be_2$. Hence, we have
$$
z_n = \eta c_1^{l - f(c_1)} \be_2 c^{kf(c)} \be_2^\ast (c_1^\ast)^l \eta^\ast = \eta \be_2 c^{l - f(c)} c^{kf(c)} (c^\ast)^l \be_2^\ast \eta^\ast = (\eta \be_2 c^l) c^{(k-1)f(c)} ((c^\ast)^l \be_2^\ast \eta^\ast)
$$
and $v(\eta \be_2 c_1^l) \subseteq W$. If $l < f(c_1)$, then we see that
$$
z_{n-1} = \al \be_2^\ast c_1^{f(c_1) - l} \eta^\ast = \al \be_3 c_1^{f(c_1) - l -1} \eta^\ast = \ga c^{kf(c)} \ga^\ast
$$
which leads to $\al \be_3 c_1^{f(c_1) - l -1} = \ga c^{kf(c)}$, $\ga = \eta$ and $c_1 = c$. So we get $\al \be_3 = \ga c^{(k-1)f(c) + l + 1}$ which means that $\al = \ga c^{(k-1)f(c) + l} \be_2$. Clearly, $c_2 = \be_3 \be_2 \in \cir c$. Hence, we have
$$
z_n = \ga c^{(k-1)f(c) + l} \be_2 \be_2^\ast (c^\ast)^l \ga^\ast = (\ga \be_2 c_2^l) c_2^{(k-1)f(c_2)} ((c_2^\ast)^l \be_2^\ast \ga^\ast)
$$
and $v(\ga \be_2 c_2^l) \subseteq W$.  \\
\textbf{Subcase 2.4.2.} $\xi = \be \xi_1$ for some $\xi_1 \in \pg$. Then we get $z_{n-1} = \al \be^\ast d^{f(d)} \be \xi_1 \eta^\ast$ and $z_n = \al \xi_1 \eta^\ast$. It follows from Lemma~\ref{Circle}(1) that there exists $c_1 \in \cir c$ such that $z_{n-1} = \al c_1^{f(c_1)} \xi_1 \eta^\ast = \ga c^{kf(c)} \ga^\ast$. So we see that $\al c_1^{f(c_1)} \xi_1 = \ga c^{kf(c)}$ and $\ga = \eta$. Noticing Lemma~\ref{CircleAdd}, one has $\al c_1^{f(c_1)} \xi_1 = \al \xi_1 c^{f(c)} = \ga c^{kf(c)}$ which leads to $\al \xi_1 = \ga c^{(k-1)f(c)}$. Therefore we obtain that $z_n = \ga c^{(k-1)f(c)} \ga^\ast$. \medskip \\
\textbf{Case 3.} $z_{n-1} = \ga (c^\ast)^{kf(c)} \ga^\ast$, where $\ga \in \pg$, $s(\ga) = v$, $v(\ga) \subseteq W$, $k \in \mathbb{Z}^+$, $r(\ga) = s(c)$ and $c \in \cnw$ with $f(c) \in \mathbb{Z}^+$. With a similar discussion as in Case~2, we can see that $z_n$ has the required form.
\end{prf}

\begin{cor} \label{Special}
The congruence $\varrho$ is special.
\end{cor}

\begin{rem}
The relation $\mathbf{R}$ used in the proof of Proposition~\ref{vinWsigma} can be chosen smaller. $\mathbf{R}$ only need to consist of all pairs $(e_v e_v^\ast, v)$ for $v \in W$ and $(c^{f(c)}, s(c))$ for $c \in C_0$ such that $f(c) \in \mathbb{Z}^+$ and $|C_0 \cap \cir d| = 1$ for every $d \in \cnw$. This can be see from (CK1).
\end{rem}

\begin{prop} \label{PCP}
$T(\varrho) = (W, f)$.
\end{prop}

\begin{prf}
Let $T(\varrho) = (W_1, f_1)$. It is not difficult to see that $W \subseteq W_1$. Conversely, if $v \in W_1$, then we certainly have $(v, e_v e_v^\ast) \in \varrho$. It follows from Lemma~\ref{vinWsigma} and the uniqueness of the expression of every nonzero element in $\In$ that $v \in W$. So we get $W_1 = W$.

For $f$ and $f_1$, arbitrarily taking $c \in \cnw$ with $f(c) \in \mathbb{Z}^+$, we observe from the definition of $\varrho$ and Lemma~\ref{Divisible} that $f_1(c) \mid f(c)$ since $(c^{f(c)}, s(c)) \in \varrho$. On the other hand, let $f_1(c) \in \mathbb{Z}^+$. Again noticing Lemma~\ref{vinWsigma} and uniqueness of the expression of every nonzero element in $\In$, we obtain that there must exist $k \in \mathbb{Z}^+$ such that $c^{f_1(c)} = c^{kf(c)}$ since $(s(c), c^{f_1(c)}) \in \varrho$. Thus, we get $f(c) \mid f_1(c)$. In conclusion, we have $f_1 = f$.
\end{prf}

\medskip
Up to now, we have completed the proof of Theorem~\ref{Main1}. To end the section, we prove Corollary~\ref{Noetherian}. \medskip

Let $(H, W, f)$ be a congruence triple. Denote the set $\{v \in V \setminus H : |s^{-1}_{G \setminus H} (v)| = 1\}$ as $\bar{V}_H$. Obviously, $W \subseteq \bar{V}_H$. Define a relation $\leq$ on the set of all congruence triples of $G$: $(H_1, W_1, f_1) \leq (H_2, W_2, f_2)$ if $H_1 \subseteq H_2$, $W_1 \setminus H_2 \subseteq W_2$ and $f_2(c) \mid f_1(c)$ for any $c \in \cnh {W_1} \cap \cnh {W_2}$. We appoint that all positive integers and $\infty$ divide $\infty$.

\begin{lem} \label{PartialOrder}
The above relation $\leq$ is a partial order on the set of all congruence triples of $G$.
\end{lem}

\begin{prf}
Obviously, the relation is reflexive and antisymmetric. To see it is transitive, we suppose that $(H_1, W_1, f_1) \leq (H_2, W_2, f_2)$ and $(H_2, W_2, f_2) \leq (H_3, W_3, f_3)$. Then we certainly have $H_1 \subseteq H_2$, $H_2 \subseteq H_3$, $W_1 \setminus H_2 \subseteq W_2$, $W_2 \setminus H_3 \subseteq W_3$, $f_2(c_1) \mid f_1(c_1)$ for any $c_1 \in \cnh {W_1} \cap \cnh {W_2}$ and $f_3(c_2) \mid f_2(c_2)$ for any $c_2 \in \cnh {W_2} \cap \cnh {W_3}$. It follows that $W_1 \setminus H_3 = (W_1 \setminus H_2) \setminus H_3 \subseteq W_2 \setminus H_3 \subseteq W_3$. If $d \in \cnh {W_1} \cap \cnh {W_3}$, then we observe that $d \in \cnh {W_2}$ because $r(d) \not\in H_3$ implies $r(d) \not\in H_2$ and $r(e) \in H_1$ implies $r(e) \in H_2$ for any exit $e$ to $d$. So we have $f_3(d) \mid f_1(d)$ which leads to $(H_1, W_1, f_1) \leq (H_3, W_3, f_3)$.
\end{prf} \medskip

By Theorem~\ref{Main2}, the mapping $T$ is a one-to-one correspondence from the set of all congruences on $\In$ to the set of all congruence triples of $G$. Concretely, given a congruence $\rho$ on $\In$, the congruence triple $T(\rho) = (H, W, f)$ is defined as what follows.
\begin{align*}
H & = 0 \rho \cap V, \\
W & = \{ s(e) \in V : (e e^\ast, s(e)) \in \rho, r(e) \not\in H\},
\end{align*}
and for any cycle $c$ such that $v(c) \subseteq W$ and every exit $e$ to $c$ satisfies $r(e) \in H$,
$$
f(c) = \left\{ \begin{array}{ll}
\infty & \mbox{if for any }k \in \mathbb{Z}^+, (c^k, s(c)) \not\in \rho, \\
m & \mbox{if } m = \min\{n \in \mathbb{Z}^+ : (c^n, s(c)) \in \rho \}.
\end{array} \right.
$$
Arbitrarily pick congruences $\rho_1, \rho_2$ on $\In$ such that $\rho_1 \subseteq \rho_2$, $T(\rho_1) = (H_1, W_1, f_1)$ and $T(\rho_2) = (H_2, W_2, f_2)$. Notice Lemma~\ref{Divisible}. It is not difficult to see that $(H_1, W_1, f_1) \leq (H_2, W_2, f_2)$. That is to say, $T$ is order preserving. Hence, if $G$ is finite, then $\In$ must be congruence Noetherian.

\section{Some other applications}

Let $S$ be a semigroup. $S$ is called {\it congruence free} if the identity relation and the universal relation are the only congruences on $S$. If $S$ has a zero element, $S$ is called {\it 0-simple} if $S$ has only $\{0\}$ and $S$ as its ideals. Let $G = (V, E, r, s)$ be a graph. A pair of vertices $u$ and $v$ in $G$ is {\it strongly connected} if there exist a path from $u$ to $v$ and a path from $v$ to $u$. A graph is {\it strongly connected} if each pair of vertices is strongly connected. Notice the correspondence of the set of ideals in $\In$ and the set of hereditary subsets of $V$. The next corollary is obvious.

\begin{cor} \label{ZeroSimple}
The following are equivalent on $\In$:
\begin{description}
\item[$(1)$] $\In$ is 0-simple;
\item[$(2)$] the only hereditary subsets of $V$ are $\emptyset$ and $V$;
\item[$(3)$] $G$ is a strongly connected graph.
\end{description}
\end{cor}

The equivalence of (1) and (3) in the following corollary is proved in \cite[Theorem~10]{MesyanMitchell}.

\begin{cor} \label{ReesCon}
The following are equivalent:
\begin{description}
\item[$(1)$] every congruence on $\In$ is a Rees congruence;
\item[$(2)$] for every hereditary subset $H$ of $V$, $\bar{V}_H = \emptyset$;
\item[$(3)$] for every $e \in E$, there exists $\al \in \pg$ such that $e$ is an exit to $\al$, $s(e) = s(\al)$ and $r(e) = r(\al)$.
\end{description}
\end{cor}

\begin{prf}
Again noticing the correspondence of the set of ideals in $\In$ and the set of hereditary subsets of $V$, we see from the definition of congruence triples, Theorem~\ref{Main2} and its proof that (1) and (2) are equivalent.

Assume that there exists some $e_0 \in E$ such that either $e_0$ is not an exit to any path, or any path $\al$ to which $e_0$ is an exit and satisfies $s(e_0) = s(\al)$ does not end at $r(e_0)$. Taking $H_0 = \{ r(\al) : \al \in \pg \setminus (V \cup \{e_0\}), s(e_0) = s(\al)\}$, we see that $H_0$ is a hereditary subset of $V$ (probably empty). However, $\bar{V}_{H_0} \neq \emptyset$ since $s(e_0) \in \bar{V}_{H_0}$. Hence, we obtain that (2) implies (3). By the definition of hereditary subsets and $\bar{V}_H$, one easily sees that (3) implies (2).
\end{prf} \medskip

The equivalence of (1) and (3) in the next corollary is proved in \cite[Theorem~3]{AshHall} for graphs where there exists at most one edge from a vertex to another and in \cite[Corollary~11]{MesyanMitchell} for general cases.

\begin{cor}
The following are equivalent:
\begin{description}
\item[$(1)$] the semigroup $\In$ is congruence free;
\item[$(2)$] $\emptyset$ and $V$ are the only hereditary subsets and $\bar{V} = \emptyset$;
\item[$(3)$] the graph $G$ is strongly connected and no vertex of $G$ has index one.
\end{description}
\end{cor}

\begin{prf}
It follows from Theorem~\ref{Main2} and Corollary~\ref{ReesCon} that (1) and (2) are equivalent. Moreover, noticing that no vertex of $G$ has index one if and only if $\bar{V} = \emptyset$, we see that (2) is equivalent to (3) from Corollary~\ref{ZeroSimple}.
\end{prf} \medskip

\noindent {\bf Acknowledgements}

I would like to express our profound gratitude to Professor Efim Zelmanov for his invaluable guidance and help.



\begin{thebibliography}{99}
\bibitem{AbramsPino} G. Abrams, G. A. Pino, The Leavitt path algebra of a graph. {\it J. Algebra} {\bf 293} (2005), 319--334.

\bibitem{AlaliGilbert} A. Alali, N. D. Gilbert, Closed inverse subsemigroups of graph inverse semigroups, {\it Comm. Algebra} {\bf 45} (11) (2017) 4667--4678.

\bibitem{AraMorenoPardo} P. Ara, M. A. Mereno, E. Pardo, Nonstable K-theory for graph algebras, {\it Algebr Represent. Theory} {\bf 10} (2007) 157--178.

\bibitem{AshHall} C. J. Ash, T. E. Hall, Inverse semigroups on graphs, {\it Semigroup Forum} {\bf 11} (1975), 140--145.

\bibitem{ChihPlessas} T. Chih, D. Plessas, Graphs ans their associated inverse semigroups, {\it Discrete Math.} {\bf 340} (10) (2017) 2408--2414.

\bibitem{Howie} J. M. Howie, {\it Fundamentals of semigroup theory}, Clarendon Press, Oxford (1995).

\bibitem{Jones} D. G. Jones, Polycyclic monoids and their generalizations, Ph.D. Thesis, Heriot-Watt University (2011).

\bibitem{JonesLawson} D. G. Jones, M. V. Lawson, Graph inverse semigroups: their characterization and completion, {\it J. Algebra} {\bf 409} (2014), 444--473.
\bibitem{Krieger} W. Krieger, On subshifts and semigroups, {\it Bull. London Math. Soc.}, {\bf 38} (2006) 617--624.

\bibitem{KumjianPaskRaeburnRenault} A. Kumjian, D. Pask, I. Raeburn, J. Renault, Graphs, groupoids, and Cuntz-Krieger algebras, {\it J. Funct. Anal.} {\bf 144} (1997) 505--541.

\bibitem{KumjianPaskRaeburn} A. Kumjian, D. Pask, I. Raeburn, Cuntz-Krieger algebras of directed graphs, {\it Pac. J. Math.} {\bf 184} (1998) 161--174.

\bibitem{MesyanMitchell} Z. Mesyan, J. D. Mitchell, The structure of a graph inverse semigroup, {\it Semigroup Forum} {\bf 93} (1) (2016) 111--130.

\bibitem{MesyanMitchellMoraynePeresse} Z. Mesyan, J. D. Mitchell, M. Morayne, Y. H. P\'{e}rese, Topological graph inverse semigroups, {\it Topology Appl.} {\bf 208} (2016) 106--126.

\bibitem{Paterson} A. L. T. Paterson, {\it Graph inverse semigroups, groupoids and their C$^\ast$-Algebras}, Birkh\"{a}user, 1999.

\bibitem{Preston} G. B. Preston, Inverse semigroup, {\it J. London Math. Soc.} (1954) 396--403.

\end{thebibliography}
\end{document}